\documentclass[]{amsart}
\usepackage{amssymb,amsmath}
\usepackage[mathscr]{euscript}

\newcounter{sec}

\newcounter{punct}[sec]

\def\punct{\refstepcounter{punct}{\arabic{sec}.\arabic{punct}.  }}

\newtheorem{theorem}{Theorem}[sec]
\newtheorem{proposition}[theorem]{Proposition}

\newtheorem{lemma}[theorem]{Lemma}

\newtheorem{corollary}[theorem]{Corollary}

\def\COUNTERS{\addtocounter{sec}{1}
              \setcounter{punct}{0}
          \setcounter{equation}{0}
          \setcounter{theorem}{0}
         }
          
          \def\sm{\smallskip}
          
\begin{document}

\newcommand{\supp}{\mathop {\mathrm {supp}}\nolimits}
\newcommand{\rk}{\mathop {\mathrm {rk}}\nolimits}
\newcommand{\Aut}{\mathop {\mathrm {Aut}}\nolimits}
\newcommand{\Out}{\mathop {\mathrm {Out}}\nolimits}
\newcommand{\OO}{\mathop {\mathrm {O}}\nolimits}
\renewcommand{\Re}{\mathop {\mathrm {Re}}\nolimits}
\newcommand{\ch}{\cosh}
\newcommand{\sh}{\sinh}

\def\0{\mathbf 0}

\def\ov{\overline}
\def\wh{\widehat}
\def\wt{\widetilde}

\renewcommand{\rk}{\mathop {\mathrm {rk}}\nolimits}
\renewcommand{\Aut}{\mathop {\mathrm {Aut}}\nolimits}
\renewcommand{\Re}{\mathop {\mathrm {Re}}\nolimits}
\renewcommand{\Im}{\mathop {\mathrm {Im}}\nolimits}
\newcommand{\sgn}{\mathop {\mathrm {sgn}}\nolimits}
\newcommand{\Isoc}{\mathop {\mathrm {Isoc}}\nolimits}
\newcommand{\PIsoc}{\mathop {\mathrm {PIsoc}}\nolimits}

\newcommand{\Sch}{\mathop {\mathrm {Sch}}\nolimits}
\newcommand{\sch}{\mathop {\mathrm {sch}}\nolimits}
\newcommand{\Fr}{\mathop {\mathrm {Fr}}\nolimits}
\newcommand{\Op}{\mathop {\mathrm {Op}}\nolimits}
\newcommand{\Mat}{\mathop {\mathrm {Mat}}\nolimits}
\newcommand{\Exp}{\mathop {\mathrm {Exp}}\nolimits}

\def\bfa{\mathbf a}
\def\bfb{\mathbf b}
\def\bfc{\mathbf c}
\def\bfd{\mathbf d}
\def\bfe{\mathbf e}
\def\bff{\mathbf f}
\def\bfg{\mathbf g}
\def\bfh{\mathbf h}
\def\bfi{\mathbf i}
\def\bfj{\mathbf j}
\def\bfk{\mathbf k}
\def\bfl{\mathbf l}
\def\bfm{\mathbf m}
\def\bfn{\mathbf n}
\def\bfo{\mathbf o}
\def\bfp{\mathbf p}
\def\bfq{\mathbf q}
\def\bfr{\mathbf r}
\def\bfs{\mathbf s}
\def\bft{\mathbf t}
\def\bfu{\mathbf u}
\def\bfv{\mathbf v}
\def\bfw{\mathbf w}
\def\bfx{\mathbf x}
\def\bfy{\mathbf y}
\def\bfz{\mathbf z}

\def\bfA{\mathbf A}
\def\bfB{\mathbf B}
\def\bfC{\mathbf C}
\def\bfD{\mathbf D}
\def\bfE{\mathbf E}
\def\bfF{\mathbf F}
\def\bfG{\mathbf G}
\def\bfH{\mathbf H}
\def\bfI{\mathbf I}
\def\bfJ{\mathbf J}
\def\bfK{\mathbf K}
\def\bfL{\mathbf L}
\def\bfM{\mathbf M}
\def\bfN{\mathbf N}
\def\bfO{\mathbf O}
\def\bfP{\mathbf P}
\def\bfQ{\mathbf Q}
\def\bfR{\mathbf R}
\def\bfS{\mathbf S}
\def\bfT{\mathbf T}
\def\bfU{\mathbf U}
\def\bfV{\mathbf V}
\def\bfW{\mathbf W}
\def\bfX{\mathbf X}
\def\bfY{\mathbf Y}
\def\bfZ{\mathbf Z}
\def\bfT{\mathbf T}

\def\frF{\mathfrak F}
\def\frD{\mathfrak D}
\def\frX{\mathfrak X}
\def\frS{\mathfrak S}
\def\frZ{\mathfrak Z}
\def\frL{\mathfrak L}
\def\frG{\mathfrak G}
\def\frg{\mathfrak g}
\def\frh{\mathfrak h}
\def\frf{\mathfrak f}
\def\frl{\mathfrak l}
\def\frp{\mathfrak p}
\def\frq{\mathfrak q}
\def\frr{\mathfrak r}
\def\frd{\mathfrak d}
\def\frb{\mathfrak b}

\def\bfw{\mathbf w}

\def\R {{\mathbb R }}
 \def\C {{\mathbb C }}
  \def\Z{{\mathbb Z}}
\def\K{{\mathbb K}}
\def\N{{\mathbb N}}
\def\Q{{\mathbb Q}}
\def\A{{\mathbb A}}
\def\U{{\mathbb U}}

\def\T{\mathbb T}
\def\P{\mathbb P}

\def\cD{\EuScript D}
\def\cL{\mathscr L}
\def\cK{\EuScript K}
\def\cM{\EuScript M}
\def\cN{\EuScript N}
\def\cP{\EuScript P}
\def\cQ{\EuScript Q}
\def\cR{\EuScript R}
\def\cT{\EuScript T}
\def\cW{\EuScript W}
\def\cY{\EuScript Y}
\def\cF{\EuScript F}
\def\cG{\EuScript G}
\def\cZ{\EuScript Z}
\def\cI{\EuScript I}
\def\cB{\EuScript B}
\def\cA{\EuScript A}
\def\cO{\EuScript O}
\def\cE{\EuScript E}
\def\ex{\EuScript E}

\def\bbA{\mathbb A}
\def\bbB{\mathbb B}
\def\bbD{\mathbb D}
\def\bbE{\mathbb E}
\def\bbF{\mathbb F}
\def\bbG{\mathbb G}
\def\bbI{\mathbb I}
\def\bbJ{\mathbb J}
\def\bbL{\mathbb L}
\def\bbM{\mathbb M}
\def\bbN{\mathbb N}
\def\bbO{\mathbb O}
\def\bbP{\mathbb P}
\def\bbQ{\mathbb Q}
\def\bbS{\mathbb S}
\def\bbT{\mathbb T}
\def\bbU{\mathbb U}
\def\bbV{\mathbb V}
\def\bbW{\mathbb W}
\def\bbX{\mathbb X}
\def\bbY{\mathbb Y}

\def\kappa{\varkappa}
\def\epsilon{\varepsilon}
\def\phi{\varphi}
\def\le{\leqslant}
\def\ge{\geqslant}

\def\B{\mathrm B}

\def\la{\langle}
\def\ra{\rangle}

\def\lambdA{{\boldsymbol{\lambda}}}
\def\alphA{{\boldsymbol{\alpha}}}
\def\betA{{\boldsymbol{\beta}}}
\def\mU{{\boldsymbol{\mu}}}

\def\const{\mathrm{const}}
\def\rem{\mathrm{rem}}
\def\even{\mathrm{even}}
\def\SO{\mathrm{SO}}
\def\SL{\mathrm{SL}}
\def\SU{\mathrm{SU}}
\def\GL{\operatorname{GL}}
\def\End{\operatorname{End}}
\def\Mor{\operatorname{Mor}}
\def\Aut{\operatorname{Aut}}
\def\inv{\operatorname{inv}}
\def\red{\operatorname{red}}
\def\Ind{\operatorname{Ind}}
\def\dom{\operatorname{dom}}
\def\im{\operatorname{im}}
\def\md{\operatorname{\,mod\,}}
\def\St{\operatorname{St}}
\def\Ob{\operatorname{Ob}}
\def\PB{{\operatorname{PB}}}
\def\Tra{\operatorname{Tra}}

\def\ZZ{\mathbb{Z}_{p^\mu}}
\def\F{\mathbb{F}}

\def\cH{\EuScript{H}}
\def\cQ{\EuScript{Q}}
\def\cL{\EuScript{L}}
\def\cX{\EuScript{X}}

\def\Di{\Diamond}
\def\di{\diamond}

\def\fin{\mathrm{fin}}
\def\ThetA{\boldsymbol {\Theta}}

\def\0{\boldsymbol{0}}

\def\FF{\,{\vphantom{F}}_3F_2}
\def\HH{\,\vphantom{H}^{\phantom{\star}}_3 H_3^\star}
\def\Ho{\,\vphantom{H}_2 H_2}

\def\disc{\mathrm{disc}}
\def\cont{\mathrm{cont}}
\def\fan{\vphantom{|^|}}

\def\osigma{\ov\sigma}
\def\ot{\ov t}

\def\Afr{\mathrm{Afr}}
\def\fr{\mathfrak{fr}}
\def\Fr{\mathrm{Fr}}

\def\tri{|\!|\!|}

\def\Diff{\mathbf{Diff}}
\def\diff{\boldsymbol{\mathfrak{diff}}}
\def\vvect{\mathfrak{vect}}
\def\vect{\boldsymbol{\mathfrak{vect}}}
\def\d{\boldsymbol{\mathfrak{vect}}}
\def\D{\mathbf {Diff}}
\def\g{\ov{\boldsymbol{\mathfrak{g}}}}

\def\G{\bfG}

\begin{center}
\bf\large


Some remarks on the group of formal diffeomorphisms of the line 


\bigskip

\sc Yury A. Neretin%
\footnote{The work is supported by the grant of FWF (Austrian Scientific Funds), PAT5335224.}
\end{center}

{\small Consider a strictly positively graded finitely generated infinite-dimensional 
real  Lie algebra $\mathfrak{g}$. It has a well-defined Lie group $\ov{\mathbf{G}}$,
which is an inverse limit of finite-dimensional nilpotent Lie groups (a pro-unipotent group). Generally, representations (even finite-dimensional representations) of $\mathfrak{g}$ and actions of $\mathfrak{g}$
on manifolds do not admit liftings to  $\ov{\mathbf{G}}$.
 There is a canonically defined dense subgroup 
$\mathbf{G}^\circ\subset \ov{\mathbf{G}}$ with a stronger (Polish) topology, which admits lifting of representations of 
$\mathfrak{g}$ in finite-dimensional spaces (and, more generally,
of representations of $\mathfrak{g}$ by bounded operators in Banach spaces).
We describe this completion for the group $\ov{\mathbf{Diff}}$ of formal
diffeomorphisms of the line, i.e., substitutions of the form 
$x\mapsto x+ p(x)$,
where $p(x)=a_2 x^2+\dots$ are formal series, and show that the group 
$\mathbf{Diff}^\circ$ consists
of series with subfactorial growth of coefficients. }

\section{Introduction}

\COUNTERS

{\bf \punct  Graded nilpotent Lie algebras and corresponding groups.%
\label{ss:1.1}} We consider
finitely generated strictly positively graded real Lie algebras $\frg$,
$$
\frg=\bigoplus_{n=1}^\infty \frg^{(n)}, \qquad \bigl[\frg^{(n)},\frg^{(m)}\bigr]
\subset \frg^{(n+m)}.
$$
Consider the completion $\g$ of $\frg$  consisting of formal series 
$\sum_{j>0} x_j$, where $x_j\in \frg^{(j)}$. 

Denote by $U(\frg)$ the universal enveloping algebra of $\frg$.
It has a natural grading generated by the grading of $\frg$, $U(\frg)=\oplus_{n\ge 0} U^{(n)}(\frg)$.
All homogeneous components $U^{(n)}(\frg)$ are finite-dimensional.
Consider the completion $\ov\bfU(\frg)$ of $U(\frg)$ consisting of formal series 
$\sum_{j\ge 0} u_j$, where $u_j\in U^{(j)}(\frg)$. It is equipped 
with the topology of component-wise convergence, 
$\sum_{j\ge 0} u_j^{(k)}$ converges to $\sum_{j\ge 0} u_j$
if for each $j$ we have the convergence $u_j^{(k)}\to u_j$
in $U(\frg)^{(j)}$.

We have a well-defined exponential map
$\g\to \ov\bfU(\frg)$: 
$$
\exp(x)=\sum_{j\ge 0}\frac{x^n}{n!} .
$$ 
By the Campbell--Hausdorff formula (see, e.g., \cite{Reu}, Ch.~3),
the subset 
$$\ov\bfG:=\exp(\g)\subset \ov\bfU(\frg)$$
 is closed with respect to multiplication,
so we get a 'Lie group' $\ov\bfG$ corresponding to the Lie algebra $\g$.

The zoo of such Lie algebras and groups includes,  in particular,
the following wide natural  families.

\sm 

{\sc A. Malcev completions of residually nilpotent discrete groups.} Consider a torsion-free  finitely generated discrete group $\Gamma$. Denote by 
$\Gamma^{[2]}:=\{\Gamma,\Gamma\}$ its commutant,  let
$\Gamma^{[3]}:=\{\Gamma^{[2]},\Gamma\}$, $\Gamma^{[4]}:=\{\Gamma^{[3]},\Gamma\}$ etc. Let some $\Gamma^{[N]}$ be trivial, i.e., let  $\Gamma$
be {\it nilpotent}.
Then, according Malcev \cite{Mal1}, $\Gamma$ admits a canonical embedding
to a connected simply connected nilpotent Lie group $G$ as a cocompact lattice.
Next, let $\Gamma$ be {\it residually nilpotent}, i.e., $\cap_N \Gamma^{[N]}$
be trivial. For each quotient $\Gamma_N:=\Gamma/\Gamma^{[N]}$
we take the corresponding Lie group $G_N$. By the functoriality of the
Malcev construction, surjective homomorphisms
$\Gamma/\Gamma^{[N+1]}\to \Gamma/\Gamma^{[N]}$ extend to surjective
homomorphisms $G_{N+1}\to G_N$. 
 According Malcev \cite{Mal2}, we
take the inverse limit
$$
\ov\bfG:=\lim\limits_\leftarrow G_N,
$$ 
and get the infinite-dimensional group $\ov\bfG$ and its Lie algebra $\g$.

The initial example of Malcev was the free group $F_n$,
in this case we come to the free Lie algebra $\frf\frr_n$ and the corresponding Lie group $\ov {\bf Fr}_n$. It
seems that the most  interesting example in this context  is the {\it Kohno
Lie algebra} $\frb\frr_n$ arising from the group of pure braids
\cite{Koh1}, \cite{Koh0}.
Recall that the  Kohno algebra has generators $L_{ij}$, where $i$, $j\le n$, $i\ne j$, $L_{ij}=L_{ji}$, and relations
$$
[L_{ij},L_{ik}+L_{jk}]=0, \qquad [L_{ij}, L_{km}]=0,
$$
for pairwise different $i$, $j$, $k$, $m$.
This algebra appears in numerous problems of representations theory
and mathematical physics. In particular, it acts in spaces of multiplicities of $n$-times tensor
 products of finite-dimensional representations (or Verma modules) of  simple Lie algebras (see, e.g., \cite{Var}), in zero weight spaces of representations
 of $\mathfrak{sl}_n$ \cite{TL}, in spaces of multiplicities of 
  quasiregular representations
 of symmetric groups $S_M$ in  functions
on  $S_M/(S_{\alpha_1}\times \dots \times S_{\alpha_n})$,
where $\sum\alpha_j=M$,
  \cite{Ner-young};
 it has very interesting geometric actions \cite{Klya}, see more references in \cite{Ner1}. 
 
 On the other hand, the Kohno algebra has numerous interesting relatives,
 see, e.g., \cite{Koh0}, \cite{Hai}, \cite{Bar}, \cite{Mar}.
 
\sm

{\sc B.  Nilpotent subalgebras of simple algebras.} Consider standard simple infinite-dimensional Lie algebras, namely Cartan
algebras of formal vector fields and affine Lie algebras 
(see \cite{Kac1}), consider their nilpotent
subalgebras. In these cases, we a priory know transparent description
of 
 groups $\ov \bfG$ (usually, this is not so). More generally, we can take nilpotent subalgebras 
of Kac--Moody algebras, see \cite{Kac2}, Ch. 1.

\sm

{\sc C.  Narrow algebras.} Infinite-dimensional nilpotent Lie algebras $\frg$ with
$\dim \frg^{(k)}\le 1$ also were  a subject of many efforts. A more general topic are algebras with uniformly bounded $\dim \frg^{(k)}$.
See, e.g., \cite{Gon}, \cite{Fia}, \cite{Koch}, \cite{ShZ}, \cite{Mil}.
 
\sm

{\bf \punct Stronger topologies.%
\label{ss:1.2}}
 A group $\ov\G$ consist of formal series. So it can act only in spaces of formal series  or actually its  actions are reduced to some
 finite-dimensional quotients of
 $\ov\G$. However, for Malcev completions there are numerous   explicit finite dimensional representations
 of algebras $\frg$ and numerous explicit actions of these algebras on  manifolds. These actions can not be integrated to $\ov \G$.
 The picture  is even more unperfect, since  monodromy operators 
 for Knizhnik--Zamolodchikov connections  are ordered exponentials
 of paths in the Kohno algebra (see, e.g., \cite{Koh2}).
  So these monodromy operators
  looks as operators of representations 
 of the Kohno Lie group.   
 
In \cite{Ner1} for algebras $\frg$ {\it generated by the
homogeneous component $\frg^{(1)}$}, there were defined two types of dense subgroups $\bfG^!\subset \bfG^\circ$
 in $\ov\G$ (with stronger topologies),  
   which allow such integrations (the smaller group $\bfG^!$ for actions of $\frg$ on compact real analytic manifolds, and the larger group $\bfG^\circ$ for finite-dimensional
   representations of $\frg$). 
   
In fact, the subgroup  $\G^\circ\subset \ov\G$
satisfies the following  property:

\sm

 --- {\it any representation $\rho$ of the Lie algebra $\frg$ in  a Banach space
by bounded operators  lifts to a representation
of $\G^\circ$.}

\sm 

We must emphasize that in standard representation theories,
representations of Lie algebras by bounded operators in infinite-dimensional Banach spaces occur very rare%
\footnote{See the basic  negative observation in \cite{Sin}.
We present a non-obvious positive example: 
 Consider the group $\cO_1$ of orthogonal operators
$g$ in a real Hilbert space  such that $g-1$ is contained in the trace class. A spinor representation of $\cO_1$ is uniformly continuous and the corresponding Lie algebra consisting of trace class skew-symmetric 
operators acts by bounded operators (this is clear from the standard
Berezin realization of the spinor representation in the fermionic Fock
space, see, e.g., \cite{Ner-book}, Chapter 4). 
This representation  can be used for a further 
production of examples by the  Olshanski technology, see
\cite{Ner-book},  Sect.~IX.5.
See a wider discussion of bounded  representations of Lie algebras in \cite{BS}.}.
 However, in our context, basic examples of representations
are finite-dimensional; on the other hand, our groups $\G^\circ$ differ from the usual 
objects of representation theories (as Lie groups, infinite-dimensional classical groups, the group of diffeomorphisms of the circle and loop groups).

\sm 

In \cite{Ner2} it was shown that the group $\mathbf{Fr}_n^\circ$
corresponding to the free Lie algebra $\mathfrak{fr}_n$
is universal in the following sense:  each connected finite-dimensional
Lie group is a quotient group of $\mathbf{Fr}^\circ_n$.

\sm 

The condition '{\it for algebras generated by $\frg^{(1)}$}' mentioned above
can be easily omitted (see Subsect. \ref{ss:last}).

The present paper is a kind of experiment: we describe the group
$\Diff^\circ$
obtained by application of this construction  to the Lie algebra of 
formal vector fields on the line
and get a strange and non-obvious topology. 

\sm

{\bf \punct The group of formal diffeomorphisms.\label{ss:1.3}}
Consider the Lie algebra%
\footnote{Its standard notation is $L_1$ or $L_1(1)$.} $\vvect$
with a basis $L_j$, where $j=1$, 2, \dots,
and relations%
\footnote{Recall that for the Virasoro algebra the commutation relations are similar,
but $j\in\Z$.}
$$
[L_n,L_m]=(m-n)L_{m+n}.
$$ 
We embed this Lie algebra to the algebra of differential operators
assuming
$$
L_n=x^{n+1}\frac{d}{dx}.
$$
The algebra $\vvect$ was a subject of many investigations, the main topics
of interest were its cohomologies (including central extensions and deformations, which 
are related to the topic of our discussion in Subsect. 1.1),
see \cite{Gon}, \cite{Fuk}, \S 2.3, \cite{FF}, \cite{Wei}.

The corresponding Lie group $\ov \G$ (see Subsect. 1.1) is the group 
$\ov \Diff$ of formal diffemorphisms of the line.
A {\it formal diffeomorphism} of $\R$ is a formal series
$$\gamma(x)=x+\sum_{j\ge 2} a_j x^j,\qquad
\text{where $a_j\in \R$.}
$$
Composition of formal diffeomorphisms
is a substitution $\gamma_1\circ \gamma_2(x)=
\gamma_1(\gamma_2(x))$ of formal series. 
The topology on $\ov \Diff$ is the usual topology on
the space of formal series, i.~e., the topology of coefficient-wise 
convergence.

This group has a canonical family of outer automorphisms
$\cE_\sigma$:
$$
\cE_\sigma \Bigl(x+\sum_{j\ge 2}a_j x^j\Bigr)=
 \Bigl(x+\sum_{j\ge 2} \sigma^{j-1}a_j x^j\Bigr).
$$
Equivalently, 
$$
\cE_\sigma\gamma(x)=\sigma^{-1} \gamma(\sigma x),
$$
i.e., this is the conjugation by the transformation
$x\mapsto \sigma x$, which is not contained in $\ov \Diff$.

Apparently, the group of formal diffeomorphisms was
introduced by Jennings \cite{Jen}, 1954, it was a subject of many works
(see, e.g., Babenko \cite{Bab} and numerous references in that survey,
see also an interesting addition to the picture in \cite{BB}).

The group $\ov \Diff$ has a natural subgroup $\Diff^\omega$
consisting of germs of analytic diffeomorphisms at 0.
An element $\gamma$ is contained in this group if
$$
\sum_{j\ge 2} |a_j|\epsilon^j<\infty \qquad
\text{for some $\epsilon>0$.}
$$
Here literature is even wider than for $\ov \Diff$ itself,
since $\Diff^\omega$ is an element of  theory of normal forms.  Again, initial references can be found in \cite{Bab}. 

Next, consider the space $W_1$ consisting of sequences
$(a_2,a_3,\dots)$ satisfying the condition
$$
\|a\|:=
\sum_{j\ge 2}\frac{|a_j|}{j!}<\infty.
$$
This space is a Banach space; in the obvious way this space is isomorphic
to the classical Banach space $\ell^1$. Denote by $\Diff_1$
the set of all formal diffeomorphisms 
$\gamma(x)=x+\sum a_j x_j^2+\dots$ such that $(a_2,a_3\dots)\in W_1$.

Denote by $\vect_1\subset \ov\vect$ the Lie algebra consisting of formal series
$\sum_{j\ge 1} b_{j} L_j$, where $(b_1,b_2,\dots)\in W_1$ (emphasize that 
the numeration of coordinates is shifted).

Our purpose is the following statements:

\begin{theorem}
\label{th:}
{\rm a)}
The set $\Diff_1$ is a subgroup in $\ov\Diff$.

\sm

{\rm b)}
 The multiplication
in $\Diff_1$ induces an analytic map $W_1\times W_1\to W_1$, and the inversion
$\gamma\mapsto \gamma^{-1}$
is an analytic map%
\footnote{This means that $\Diff_1$ is a Banach--Lie group.
Many elementary statements about finite-dimensional Lie groups 
survive for  Banach-Lie groups, such groups were a subject of optimism
in the 1960-s, and even Nikola Bourbaki \cite{Bou} started his exposition
of theory of Lie groups in this generality. Quite soon, it became clear that 'infinite-dimensional Lie groups' arising in different domains of mathematics and in mathematical physics usually are not Banach-Lie groups. 'Infinite-dimensional Lie groups' are 'manifolds' modeled by locally convex spaces of more general types (as Fr\'echet spaces or inductive limits of Fr\'echet spaces) or are 'manifolds' only in heuristic sense (as the complete unitary group
of a Hilbert space). However, the property 'a group
is  Banach--Lie group' is a strong positive property, and constructions 
 of \cite{Ner1}, \cite{Ner2} with the present paper enlarge a zoo of Banach--Lie groups.}$^,$%
 \footnote{Generally, formal series with such growth of coefficients  diverge everywhere, but such series are usual in asymptotic analysis.}
  $W_1\to W_1$.
\end{theorem}

\begin{proposition}
\label{pr:exp}
{\rm a)} There is a well-defined analytic exponential map
$\vect_1\to \Diff_1$.

\sm 

{\rm b)} There is a well-defined analytic
 logarithmic map from a neighborhood
of unity in $\Diff_1$ to  $\vect_1$.
\end{proposition}
  
We can slightly vary the construction considering the space
$W_\sigma$, where $\sigma>0$ is a parameter, consisting of sequences
$a=(a_2,a_3,\dots)$ satisfying the condition
$$
\|a\|_\sigma:=\sum_{j\ge 2}\frac{a_{j}\sigma^{j-1}}{j!}<\infty,
$$
 the corresponding group $\Diff_\sigma$, and the corresponding Lie
 algebra $\vect_\sigma$.
  The outer automorphism
$\cE_\sigma$ identify $\Diff_1$ with $\Diff_\sigma$
and $\vect_1$ with $\vect_\sigma$.

We also define the group
$$
\Diff^\circ:=\cap_{\sigma} \Diff_\sigma
$$
and the Lie algebra
$$
\vect^\circ=\cap_{\sigma} \vect_\sigma.
$$

The previous proposition implies:

\begin{proposition}
There is a well-defined exponential map $\vect^\circ\to \Diff^\circ$.
\end{proposition}

Also, we show that:

\begin{proposition}
\label{pr:universality}
Any representation of the Lie algebra $\vvect$ by bounded operators
in a Banach space can be canonically integrated to a uniformly continuous representation 
of $\Diff^\circ$ in the same space.
\end{proposition}

{\bf\punct The further structure of the paper.}
In Section \ref{s:2}, we verify Theorem \ref{th:} and Proposition
\ref{pr:exp}.
In Section \ref{s:3}, we show that $\Diff^\circ$
is a special case of  a construction
of \cite{Ner1} and verify Proposition \ref{pr:universality}.

\section{The group $\Diff_1$\label{s:2}}

\COUNTERS



{\bf \punct The space $\ov V$ of formal series.}
We consider the linear space $\ov V$ consisting of all formal series 
$$
F(x)=\sum_{m=0}^\infty u_m x^m
$$
 with real coefficients $u_m$. 

We say that an operator $S$ is {\it triangular} (resp., {\it strictly  triangular}) 
if its matrix elements $s_{ij}$ in the basis $x^k$ are zeros for
$i<j$ (resp., $i\le j$).

Equivalently,  consider the following filtration
$$\ov V\supset x\, \ov V\supset x^2\,\ov V\supset\dots $$
An operator $S$ is triangular if 
$S( x^k\, \ov V) \subset  x^k\, \ov V $
 for all $k$ and strictly triangular if
 $S( x^k\, \ov V) \subset  x^{k+1}\, \ov V $.
 So a triangular operator $S:\ov V\to\ov V$ determines a sequence
 of linear transformations 
 $$ \ov V/x^k \,\ov V \to \ov V/x^k\, \ov V, $$
 and $S$ is determined by this sequence.

\sm 
 
 For any formal diffeomorphism 
 $\gamma=x+a_2 x^2+\dots$,
 we define the linear  operator 
 $T(\gamma):
 \ov V\to \ov V$ by
$$
T(\gamma)\, F(x)=T(\gamma)\,\Bigl(\sum_{m\ge 0} u_m x^m\Bigr):=
\sum_{m\ge 0} u_m (x+a_2 x^2+\dots)^m.
$$
Obviously, we get a representation of $\ov\Diff$:
$$
T(\gamma_1)\, T(\gamma_2)=
T(\gamma_1\circ \gamma_2).
$$
Clearly, each matrix element of 
$T(\gamma)$ in the basis
$x^l$ is a polynomial in the variables $a_j$. Notice also that
the operator $T(\gamma)-1$  is strictly triangular. 

Next, consider the Lie algebra $\ov\vect$ consisting of all formal series
$$
\cL= \sum_{j\ge 1} p_j L_j= \sum_{j\ge 1} p_j x^{j+1} \frac{d}{dx}.
$$
Such operators are strictly  triangular.

Notice that for a strictly  triangular matrix $S$
we have well-defined 
$$
\exp S:=1+\sum_{k> 0}\frac{1}{k!}\, S^k, \qquad
\ln (1+S):=\sum_{k>0}\frac{(-1)^{k-1}}{k} S^k,
$$ 
and these operations are inverse one to another (indeed, this is valid
if we descent to spaces $\ov V/x^N \ov V)$.

We need the following statement \cite{Jen}:

\begin{proposition}
{\rm a)} For any $\cL\in \ov\vect$, 
the operator $\exp (\cL)$ has a form 
$T(\gamma)$ with $\gamma\in\ov \Diff$.

\sm

{\rm b)} For any $\gamma\in \ov\Diff$,
we have $\ln T(\gamma):=\ln\bigl(1+(T(\gamma)-1)\bigr)\in\ov\vect$. 
\end{proposition}




\begin{proposition}
The group $\exp(\ov\vect)\subset \ov\bfU(\vvect)$
defined in Subsect. {\rm \ref{ss:1.1}} is isomorphic to $\ov\Diff$.
\end{proposition}

{\sc Proof.} The representation of $\vvect$ in $\ov V$
determines a representation of the completed enveloping algebra
$\ov\bfU(\vvect)$
in the same space. So we have a representation of the group
$\exp(\ov\vect)$ in $\ov V$. On the other hand,
our representation of $\ov\bfU(\vvect)$ sends exponentials to exponentials.
\hfill $\square$

\sm

{\bf \punct The Lie algebra $\vect_{1}$.}
By $V_1$ we denote the subspace in $V$ consisting of formal series
$\sum_{m} u_mx^m $ satisfying
 \begin{equation}
 \sum_{m}\frac{|u_m|}{m!}<\infty .
 \end{equation}
 
This space is isometric to $\ell^1$. So, for any
linear  operator $A$ in $V$ we have (see, e.g., \cite{KA}, V.2.8.b)
$$
\|A\|=\sup_{m} \frac{\|A x^m\|}{\|x^m\|}=
\sup_{m}  m! \, \|A x^m\|.
$$

Denote by $\cB(V_1)$ the Banach algebra of all bounded operators in 
$V_1$.

\sm 


By $\vect_{1}\subset \ov\vect$ we denote
the subspace consisting of all operators 
$$
\cL=\bigl(\sum_{j\ge 1} p_j x^{j+1}\bigr)\frac{d}{dx}
$$
satisfying the condition 
\begin{equation}
\sum_j \frac{|p_j|}{(j+1)!}<\infty.
\label{eq:L-Lie}
\end{equation}

\begin{proposition}
\label{pr:lie-algebra}
An operator $\cL\in \ov\vect$ is bounded 
 in $V_1$ if and only if $L\in \vect_{1}$.
\end{proposition}

{\sc Proof.} Applying $\cL$ to a monomial $x\in V_1$,
we get
$$\|\cL x\|=\sum_j \frac{|p_j|}{(j+1)!},$$
so  condition \eqref{eq:L-Lie}  is necessary. 

Next,
let \eqref{eq:L-Lie} be satisfied. Then
\begin{multline*}
\|\cL\|
=\sup_{m\ge 0} m!\, \|\cL x^m\|=\\=
\sup_{m>0}\Bigl( m!\,\Bigl\|\sum_{j>0} p_j x^{j+m}\Bigr\|\cdot m\Bigr)
=
\sup_{m>0}\Bigl(m\cdot m! \cdot\sum_{j> 0} \frac{|p_j|}{(m+j)!} \Bigr)=
\\
=
\sup_{m>0}\Bigl(\sum_j \frac{|p_j|}{(j+1)!}\cdot
\frac{m\cdot m! }{(j+2)\dots (j+m)} \Bigr)
\le \\\le
\sum_j \frac{|p_j|}{(j+1)!}
\,\cdot\,
\sup_{m>0}\sup_{j>0} \frac{m\cdot m! }{(j+2)\dots (j+m)}
=\\=
\sum_j \frac{|p_j|}{(j+1)!}
\cdot \sup_{m>0}\, \frac{m\cdot m! }{\frac12\, (1+m)!}
=2\cdot \sum_j \frac{|p_j|}{(j+1)!}.
\qquad\square
\end{multline*}

\begin{corollary}
The $\vect_1$ is a Lie algebra.
\end{corollary}

Indeed, if $\cL_1$, $\cL_2\in \cB(V_1)$, then 
$[\cL_1, \cL_2]\in \cB(V_1)$.

\sm

{\bf \punct The group $\Diff_{1}$.}
We consider a Banach space $W_1$
consisting of sequences $(a_2,a_3,\dots)$
such that
\begin{equation}
[|a|]:=
\sum_{j\ge 2}\frac{|a_j|}{j!}<\infty
\label{eq:sum-a}
\end{equation}
and
 the subset  $\Diff_1\subset\ov \Diff$ 
consisting of $\gamma=x+a_2 x^2+\dots\in\ov\Diff$ 
such that $(a_2,a_3,\dots)\in W_1$.

\begin{theorem}
\label{th:bounded}
An operator $T(\gamma)$ is bounded in $V_1$ if and only if 
$\gamma\in \Diff_1$. The map $\gamma\mapsto T(\gamma)$ is analytic
as a map from $W_1$
to  $\cB(V_1)$.
\end{theorem}

{\sc Proof.}
Since $T(\gamma)\,x=x+\sum_{j>1}a_j x^j$,
the condition \eqref{eq:sum-a} is necessary.

Let $\gamma\in \ov \Diff$ be arbitrary. Denote $\gamma(x)-x:=x^2h(x)$.
Then we have
\begin{equation}
T(\gamma)F(x)=T\bigl(x+x^2h(x)\bigr)F(x)=
\sum_{n\ge 0}\frac {1}{n!} \bigl(x^2 h(x)\bigr)^n F^{(n)}(0).
\label{eq:T-taylor}
\end{equation}
Indeed,
if $F$ and $\gamma$ are analytic at zero, then this identity is the Taylor
formula. Since actually this is an identity for polynomials expressions,
it is valid for formal series. 

Next, we define operators
$Q_n:\ov V\to \ov V$ assuming
\begin{align*}
Q_n x^n&:=x^{2n}\frac{d^n}{dx^n}x^n=n!\, x^{2n}
\\
Q_n x^k&:=0,\quad \text{if $k\ne n$.} 
\end{align*}
We also define the operator $H:\ov V\to \ov V$ assuming
\begin{align*}
H x^k:= h(x)\,x^k, \qquad \text{if $k\ge 2$;}
\\
H\, 1=0,\qquad H\, x=0.
\end{align*}
We rewrite \eqref{eq:T-taylor} in the form
\begin{equation}
T(\gamma)=\sum_{n\ge 0} \frac{1}{n!} H^n Q_n.
\label{eq:taylor-2}
\end{equation}
Next, let $\gamma\in \Diff_1$.
Let us examine norms of $H$ and $Q_n$ in $V$.
Clearly, 
$$
\|Q_n\|=n!
\, \|Q_n x^n\|=\frac{n! \,n!}{(2n)!}<\frac{2\sqrt n}{2^{2n}}<1.
$$
Also,
\begin{multline*}
\|H\|=\sup_{k\ge 2} k! \,\|H x^k\|=\sup _{k\ge 2} 
 k! \, \|\sum_{j\ge 2} a_j x^{j+k-2}\|
= \sup _{k\ge 2} 
 k! \, \sum_{j\ge 2} \frac{|a_j|}{(j+k-2)!}
 =\\=\sup _{k\ge 2}  \sum_{j\ge 2} \frac{|a_j|}{j!}
 \frac{k!}{(j+1)\dots (j+k-2)}.
\end{multline*}
But
$$ \frac{k!}{(j+1)\dots (j+k-2)}
=1\cdot 2\cdot\frac{3}{j+1}\cdot \frac{4}{j+2}\cdot
\dots\vphantom{.}\cdot
\frac{k}{j+k-2}\le 1\cdot 2\cdot 1\cdot\vphantom{.} \dots \vphantom{.} \cdot1.
$$
Therefore, 
$$\|H\|\le 2\cdot \sum_{j\ge 2} \frac{|a_j|}{j!}<\infty.
$$
For this reason, the series \eqref{eq:taylor-2}  converges with respect to operator norm. The operator $H$ analytically depends on $\gamma$
(since $\gamma\to H$ is a linear map $W_1\to\cB(V_1)$), so
$T(\gamma)$ also analytically depends on $\gamma$.
\hfill
$\square$

\sm

\begin{corollary}
\label{cor:mult}
The subset $\Diff_1$ is closed with respect to the multiplication. 
The multiplication in $\Diff_1$ is an analytic operation.
\end{corollary}

Indeed, a product of bounded operators is a bounded operator.

\begin{corollary}
If a formal vector field $\cL$ is contained in $\vect_1$,
then its exponential is contained in $\Diff_1$.
\end{corollary}

Indeed, $\cL$ determines a bounded operator, the exponential of
a bounded operator $\cL$ is a bounded operator $T(\gamma)$, and therefore
$\gamma\in \Diff_1$.

\begin{corollary}
If $\|T(\gamma)-1\|\le 1$, then $\ln T(\gamma)\in \vect_1$.
\end{corollary}

Indeed, a logarithm of a bounded operator is a bounded operator.

\begin{theorem}
The set $\Diff_1$ is a subgroup in $\ov \Diff$.
\end{theorem}

Indeed, by Corollary \ref{cor:mult}, $\Diff_1$ is closed with respect 
to the composition. It remains to prove that 
$\Diff_1$ is closed with respect to the inversion, so we must establish
the following Lemma \ref{l:inversion}.

\sm 

{\bf \punct Inverse diffeomorphisms.} 

\begin{lemma}
\label{l:inversion}
Let $\gamma \in \Diff_1$. Then the inverse formal diffeomorphism
also is contained in $\Diff_1$.
\end{lemma}

{\sc Proof.} Let $\gamma(x)=x-a_2 x^2-a_3 x^3-\dots\in \Diff_1$.
 Let $\mu(x)=\gamma^{-1}(x)$
be the inverse formal diffeomorphism,
$$
\mu(x)=x+c_2x^2+c^3 x^3+\dots
$$
By the Lagrange inversion formula (see, e.~g., \cite{AAR}, Addendum F),
 we have
$$
c_n=\frac 1{n!}\,\Bigl[\frac{d^{n-1}}{dx^{n-1}}
\Bigl(\frac{x}{\gamma(x)}\Bigr)^n\Bigr]\Bigr|_{x=0}.
$$
Expanding $(x/\gamma(x))^n$
in a series we get
\begin{multline*}
\Bigl(\frac{x}{\gamma(x)}\Bigr)^n=
(1-a_2 x-a_3x^2- \dots)^{-n}=\\=\sum_{p\ge0}
\frac{ n(n+1)\dots (n+p-1)}{p!}\bigl(\sum_{j\ge 1} a_{j+1}\,x^j\bigr)^p
=\\=
\sum_{p\ge 0} n(n+1)\dots (n+p-1)
\sum_{\footnotesize\begin{matrix}k_1\ge 0, k_2\ge 0,
\vphantom{.}\dots\vphantom{.}:\\ \sum k_j=p\end{matrix}}\,\,\, \prod_{j\ge 1}\frac{(a_{j+1}\,x^j)^{k_j}}{k_j!}.
\end{multline*}
Then $c_n$ is $\frac 1n \sigma_{n-1}$, where $\sigma_{n-1}$
is the coefficient at $x^{n-1}$ in the last expression.
So
$$
c_n=\frac 1n \sum_{p\ge 0} n(n+1)\dots (n+p-1)
\sum_{\footnotesize\begin{matrix}k_1\ge 0, k_2\ge 0,
\vphantom{.}\dots\vphantom{.}:\\ \sum k_j=p,\, \sum jk_j=n-1\end{matrix}} \prod_{j\ge 1}\,\frac{(a_{j+1})^{k_j}}{k_j!}.
$$
We must show that 
$$
S:=\sum_{n\ge 2} \frac{|c_{n}|}{j!}<\infty.
$$
 Using the inequality $|\sum \dots|\le \sum|\dots|$
and passing to the summation in $k_1$, $k_2$, \dots
(so $p=\sum k_j$, $n=\sum jk_j+1$), we get
$$
S\le \wt S:=\sum_{k_1,k_2,\dots}
\frac{(\sum jk_j+2)\dots(\sum jk_j+\sum k_j)}{(\sum jk_j+1)!}
\cdot \prod\frac{|a_{j+1}|^{k_j}}{k_j!}.
$$
Further, we transform the latter expression to the form
\begin{equation}
\wt S=\sum_{k_1,k_2,\dots} U(k_1,k_2,\dots)
 \prod\frac 1 {k_j!} \Bigl(\frac{|a_{j+1}|}{(j+1)!}\Bigr)^{k_j},
 \label{eq:wt-S}
\end{equation}
where  
$$
U(k_1,k_2,\dots):=\frac{(\sum jk_j+2)\dots(\sum jk_j+\sum k_j)
\prod_{j}\bigl[(j+1)!\bigr]^{k_j}
}{(\sum jk_j+1)!}.
$$

\begin{lemma}
\label{l:U}
There exist constants $L$ and $M$ such that
$$
U(k_1,k_2,\dots)\le L\cdot M^{\sum k_j}
$$
for all collections $k_1$, $k_2$, \dots
\end{lemma}

{\sc Proof of Lemma \ref{l:U}.}
We transform $U(k_1,k_2,\dots)$
as
\begin{equation}
\biggl[
\frac{(\sum jk_j+2)\dots(\sum jk_j+\sum k_j)}
{(\sum (j-1)k_j+3)\dots(\sum jk_j+1)}\biggr]\,
\cdot\,\biggl\{
\frac{\prod_{j}\bigl[(j+1)!\bigr]^{k_j}}{
(\sum (j-1)k_j+2)!}\biggr\}.
\label{eq:long-product}
\end{equation}

Let us examine the first factor. Consider a fraction
$$
\Phi_l(t)=\frac{t(t+1)\dots(t+l-1)}{(t-l)(t-l+1)\dots(t-1)}.
$$
Clearly, it  decreases on the ray $t\ge l+1$, and therefore its maximal
value on this ray is
$$
\Phi_l(l+1)=\frac{(2l)!}{l!\,l!}\le \frac{2^{2l}}{\sqrt{\pi l}}.
$$
We apply this remark to the first factor in \eqref{eq:long-product},
$t=\sum jk_j+2$, $l=\sum k_j-1$ and get the upper estimate
$$
\Bigl[\dots\Bigr] \le \mathrm{const}\cdot 4^{\sum k_j}.
$$

Next, let us examine the second factor in \eqref{eq:long-product}. We write it in the form
$$
\biggl\{\dots\biggr\}=
\frac{(1\cdot 2)^{\sum_{j\ge 1} k_j}}{1\cdot 2}
\cdot \biggl( \frac {\prod_{j\ge 2} \bigl(3\cdot\dots \vphantom{.}\cdot (j+1)\bigr)^{k_j} }
{3\cdot 4\cdot \dots \vphantom{.} \cdot (\sum (j-1)k_j+2)} 
\biggr)
.$$
We claim that the factor  $\bigl(\dots\bigr)$ is $\le 1$.
 Indeed, the numbers of factors
in the numerator and the denominator coincide; namely, they are
 equal $\sum (j-1)k_j$.  In the numerator,
we write $\sum_{j>1} k_j$ products of type 
$3\cdot\vphantom{.} \dots\vphantom{.}\cdot (j+1)$ one after another,
\begin{equation}
\underbrace{(3)\cdot \dots \vphantom{.} \cdot (3)}_{\text{$k_2$ times}}
\cdot \underbrace{(3\cdot 4)\cdot \dots \vphantom{.} \cdot (3\cdot 4)}_{\text{$k_3$ times}} \cdot \underbrace{(3\cdot 4\cdot 5)\cdot \dots \vphantom{.} \cdot (3\cdot 4\cdot 5)}_{\text{$k_4$ times}}\cdot \dots\, .
\label{eq:345}
\end{equation}
Then each factor of numerator is $\le$ the corresponding factor of denominator. If  $\sum_{j>1} k_j=1$ (i.e., if a product \eqref{eq:345}
consists of one subproduct $3\cdot \dots \vphantom{.}\cdot (\alpha+1)$), then the numerator equals  the denominator. Otherwise, we have a strict inequality.

 This completes the proof of Lemma \ref{l:U}.
\hfill $\square$

\sm

{\sc End of proof of Lemma \ref{l:inversion}.}
Now we are ready to estimate \eqref{eq:wt-S}:
$$
\wt S\le  L\cdot\sum_{k_1,k_2,\dots} M^{\sum k_j}\prod_{j\ge 1}
\frac{1}{k_j!}\Bigl(\frac{|a_{j+1}|}{(j+1)!} \Bigr)^{k_j}
=L\cdot \exp\Bigl\{M\sum_{j\ge 1}\frac{|a_{j+1}|}{(j+1)!}  \Bigr\}
<\infty.\qquad \square$$

\sm

Thus, $\Diff_1$ is a group. Since the inversion of bounded operators
is an analytic operation, the inversion in $\Diff_1$ also is analytic.

\section{Comparison with the general construction\label{s:3}}

\COUNTERS

In Subsect. \ref{ss:1.3}, we defined a family
of  subgroups $\Diff_t\subset \ov\Diff$ with Banach--Lie structure,
 namely $\Diff_t=\cE_{t^{-1}}(\Diff_1)$,
$$
t>s\quad \Rightarrow \quad \Diff_t\subset\Diff_s.
$$
The general construction mentioned in 
Subsect. \ref{ss:1.2} also produces a family of 
subgroups $\Diff_{[t]}\subset \ov\Diff$ with Banach--Lie structure,
they also satisfy
$\Diff_{[t]}=\cE_{t^{-1}}(\Diff_{[1]})$.
The author does not know precise relations between these two scales;
in this section we show that they are  very close one to another,
see
\eqref{eq:sub-sub}. Also,
$$
\bigcap_{t>0}\Diff_{[t]}= \bigcap_{t>0}\Diff_t.
$$

{\bf \punct Submultiplicative norms on the free associative algebra.%
\label{ss:sub}}
Consider the free associative algebra $\cA_2$ with two generators 
$\omega_1$ and $\omega_2$. We denote monomials
as
$$\omega^I:=\omega_{i_1}\omega_{i_2}\dots\omega_{i_m}$$ 
By $d_1(I)$ (resp. $d_2(I)$) we denote the number
of entries  of $\omega_1$ (resp. $\omega_2$) to a monomial
$\omega^I$.


Fix $t_1$, $t_2>0$. We define the norm
$R_{[t_1,t_2]}$ on $\cA_2$ by
$$
R_{[t_1,t_2]}\bigl(\sum_I c_I \omega^I\bigr):=
\sum |c_I|\, t_1^{d_1(I)} t_2^{d_2(I)}.
$$ 
This norm is {\it submultplivative}, i.e.,
$$
R_{[t_1,t_2]}(u\cdot v)\le R_{[t_1,t_2]}(u)\,R_{[t_1,t_2]}(v).
$$

\begin{lemma}
{\rm a)} Let $Q$ be a submultiplicative norm
on $\cA_2$ such that $Q(\omega_1)=t_1$, $Q(\omega_2)=t_2$.
Then
$Q(u)\le R_{[t_1,t_2]}(u)$ for all $u\in\cA_2$.

\sm

{\rm b)} If $s_1\le t_1$, $s_2\le t_2$,
then $ R_{[s_1,s_2]}(u)\le  R_{[t_1,t_2]}(u)$.

\sm

{\rm c)} Let $\rho$ be a representation of $\cA_2$ 
in a Banach space $H$. Let $\|\rho(\omega_1)\|=t_1$,
 $\|\rho(\omega_2)\|=t_2$. Then for any
 $u\in\cA_2$ we have $\|\rho(u)\|\le R_{[t_1,t_2]}(u)$.
\end{lemma}

The statements a), b) are obvious. It remains to notice that
 $\|\rho(u)\|$
is a submultiplicative norm on $\cA_2$.

\sm

{\bf \punct Completions of the universal enveloping algebra  $U(\vvect)$.%
\label{ss:completions}}
 We apply the considerations of 
Subsect. \ref{ss:1.1} to $\frg=\vvect$.
Consider the universal enveloping algebra $U(\vvect)$ and its completion 
$\ov \bfU(\vvect)$  consisting of formal series
$u=\sum_{n\ge 0} u^{(n)}$, where $u^{(n)}\in U^{(n)}(\vvect)$
are homogeneous elements elements.
 For $t>0$ consider the automorphism $\cE_t$ of $U(\vvect)$
 (and of $\ov \bfU(\vvect)$),
defined by
 $$
 \cE_t \bigl(\sum_{n\ge 0} u^{(n)}\bigr):=
  \sum_{n\ge 0} t^n u^{(n)}.
 $$

Consider the homomorphism
$\pi:\cA_2\to U(\vvect)$ sending
$$
\omega_1\mapsto L_1,\qquad \omega_2\mapsto L_2.
$$
 So, $U(\vvect)$ becomes a quotient of $\cA_2$ by a certain  ideal
 (which is a sum of its homogeneous components). We equip $U(\vvect)$ 
 with a factor-norm $Q_{[1]}$
 defined by
\begin{equation}
Q_{[1]}(u):=\min_{w\in \cA_2:\,\pi (w)=u} R_{[1,1]}(w).
\label{eq:Q-min}
\end{equation}

Notice that
$$
Q_{[1]}(L_1)=1,\qquad Q_{[1]}(L_2)=1.
$$

\begin{lemma}
{\rm a)} If $u\in U^{(k)}(\vvect)$, then 
it is sufficient to take a minimum over $w\in\cA_2$
of degree $k$.

\sm

{\rm b)} $Q_{[1]}(\sum u^{(n)})=\sum Q_{[1]}( u^{(n)})$.

\sm

{\rm c)}
For any submultiplicative norm $P$ on $U(\vvect)$
such that $P(L_1)=1$, $P(L_2)=1$ we have
$P(u)\le Q_{[1]}(u)$. 
\end{lemma}

{\sc Proof.}
b) Decompose $w$ as a sum of its homogeneous components,
$w=\sum w^{(n)}$.
  The  condition 
$\pi(\sum w^{(n)})=\sum u^{(n)}$ is equivalent to
the collection of independent conditions
$\pi (w^{(n)})=u^{(n)}$ for all $n$.
So the minimum is the sum of minima.

\sm

a) is a special case of b).

\sm

c) Indeed, $P(\pi (v))$ is a submultiplicative
norm on $\cA_2$. So $P(\pi(v))\le R_{[1,1]}(v)$.
So for any  $w\in\cA_2$ such that $\pi(w)=u$
we have $P(\pi (w))\le R_{[1,1]}(w)$. By 
\ref{eq:Q-min}, we get our statement.
\hfill $\square$

\sm
 
 
For each $t>0$ we define a norm $Q_{[t]}$ on $U(\vvect)$
 assuming
 $$
 Q_{[t]}(u):=Q_{[1]}(\cE_{[t]}u).
 $$  
Consider the completion $\bfU_{[t]}(\vvect)$ of $U(\vvect)$
with respect to a norm $Q_{[t]}$, so it is a space of formal series
$\sum u^{(n)}\in \ov\bfU(\vect)$ such that $\sum t^n Q_{[1]}(u^{(n)})\le 1$.
Denote
$$
\vect_{[t]}:=\vect\cap \bfU_{[t]}(\vvect),
\qquad \Diff_{[t]}:=\exp(\ov\vect)\cap \bfU_{[t]}(\vvect).
$$
Clearly, $\exp(\vect_{[t]})\subset \Diff_{[t]}$.
Also, $\ln (\cdot)$ sends a neighborhood of
the unit in $\Diff_{[t]}$ to a neighborhood of 
zero in $\vect_{[t]}$. 

\sm

So, $\Diff_{[t]}$ is a  Banach--Lie group.

\begin{proposition}
The groups $\Diff_{[t]}$  and $\cap_t\Diff_{[t]}$ are connected
and contractible.
\end{proposition}

{\sc Proof.}
Let $\gamma\in \Diff_{[t]}$.
let $\gamma=\exp(\cL)$ (apparently, generally, $\cL\notin \vect _{[t]}$).
 The curve $\cE_\sigma(\gamma)$, where $\sigma\in[0,1]$,
connects $\gamma$ with the unit. By definition, 
$$\cE_\sigma(\gamma)\in \Diff_{[\sigma t]}\subset \Diff_{[t]}\subset
\bfU_{[t]}(\vvect).$$
On the other hand,
$$
\cE_\sigma(\gamma)=\exp\bigl(\cE_\sigma (\cL)\bigr)
\in \exp(\vect).
$$
Therefore, $\cE_\sigma(\gamma)\in \Diff_{[t]}$.
Also, we get a homotopy $[0,1]\times \Diff_{[t]}\to \Diff_{[t]}$
defined by $(\sigma,\gamma)\to \cE_\sigma (\gamma)$.
\hfill $\square$

\sm


{\bf \punct Representations of $\Diff_{[t]}$.%
\label{ss:representations}}
Consider a representation $\rho$  of the Lie algebra
$\vvect$ in a Banach space $H$ such that norms
of the generators satisfy 
$$
\tri \rho(L_1)\tri\le t,\qquad \tri \rho(L_2)\tri\le t^2.
$$
Then the algebra $\bfU_{[t]}(\vvect)$ acts by bounded operators. The operator
norm  $\tri\cdot\tri$ induces a submultiplicative norm
$\tri\rho(u)\tri$ 
on $U_{[t]}$ and therefore
\begin{equation}
\tri \rho(u)\tri\le Q_{[t]}(u).
\label{eq:sub}
\end{equation}
In particular, we get a representation 
of the Lie algebra $\vect_{[t]}$ and a representation 
of the group $\Diff_{[t]}$ in $H$ (since both objects are contained in
$\bfU_{[t]}(\vvect)$).

\sm 

{\bf\punct Upper estimates of norms  $Q_{[t]}$ of elements of $\vvect$.}

\begin{lemma}
For any $\cL=\sum p_j L_j\in \vvect$, we have
\begin{equation}
Q_{[t]}\bigl(\sum p_j L_j\bigr)\le |p_1|+
\frac 14\sum_{j>1} \frac{|p_j| (2t)^j}{(j-2)!}.
\label{eq:Q}
\end{equation}
\end{lemma}

{\sc Proof.}  
We have
$$
L_{n}=\frac1{n-2}[L_1,L_{n-1}].
$$
Since $Q_{[t]}(L_1)=t$, 
$$
Q_{[t]}(L_n)\le  \frac1{n-2} Q_{[t]}(L_1 L_{n-1}-L_{n-1}L_1)
\le \frac{2t}{n-2} Q_{[t]}(L_{n-1}).
$$
Iterating, we get
$$
Q_{[t]}(L_n)\le \frac{2^{n-2}t^n}{(n-2)!},
$$
and this implies our statement.
\hfill $\square$

\begin{corollary}
For any $\epsilon>0$ we have
$\vect_{[t]}\supset \vect_{2t+\epsilon}$.
\end{corollary}

{\sc Proof.} Denote by $\Xi(\cL)$ the right hand side of \eqref{eq:Q}.
Let $\cL=\sum p_j L_j\in \vect_{2t+\epsilon}$.
Then the following sum $\Sigma(\cL)$ is finite:
\begin{multline*}
\Sigma(\cL)=
\sum_{j\ge 1}\frac{|p_j|(2t+\epsilon)^j}{(j+1)!}
=\\=
\frac{|p_1|(2t+\epsilon)}{2}+\sum_{j\ge 2}\frac{|p_j|(2t)^j}{(j-2)!}\cdot
\biggl\{\frac{1}{(j-1)j(j+1)} 
\Bigl(\frac{2t+\epsilon}{2t}\Bigr)^j\biggr\}.
\end{multline*}
The expression $\{\dots\}$ tends to $\infty$ as $j\to\infty$.
 Therefore $\min_j\{\dots\}>0$.
 Hence, there exists $\delta=\delta(t,\epsilon)>0$ such that
 $\Sigma(\cL)\ge \delta\, \Xi(\cL)$. But $\Xi(\cL)\ge Q_{[t]}(\cL)$.
 \hfill $\square$

\sm

{\bf\punct Lower estimates of norms  $Q_{[t]}$ of elements of $\vvect$.}

\begin{lemma}
\label{cor:1}
For any
$\sum p_j L_j\in \vvect$,
we have
$$Q_{[t/2]}(\sum p_j L_j)\ge \sum \frac{|p_j| t^j}{(j+1)!}.$$
\end{lemma}

{\sc Proof.} Consider the Banach space $V_t$
consisting of formal series $F(x)=\sum_{m\ge 0} u_m x^m$ such that
$$
\|F\|_t:=\sum \frac{|u_m|\, t^m}{m!}<\infty.
$$
 As in the proof of Proposition \ref{pr:lie-algebra},
 $$
\sum \frac{|p_j|t^{j}}{(j+1)!}\le  \bigl\|\sum p_j L_j\bigr\|_t\le
2 \sum \frac{|p_j|t^{j}}{(j+1)!}.
 $$

Also, we have
$$
\|L_1\|_t=\frac12 t,\qquad \|L_2\|_t=\frac 16 t^2< 
\Bigl(\frac12 t\Bigr)^2;
$$
 Now our statement follows from \eqref{eq:sub}.
\hfill $\square$

\begin{corollary}
\label{cor:2}
$\vect_{[s]}\subset \vect_{2s}$.
\end{corollary}

Thus,
\begin{equation*}
\vect_{2t}\supset\vect_{[t]}\supset \vect_{2t+\epsilon}.
\end{equation*}
Since we have exponential maps
$\vect_{[s]}\to \Diff_{[s]}$, $\vect_s\to \Diff_s$
and groups $\Diff_{[s]}$, $\Diff_s$ are connected,
we get similar inclusions for the groups
\begin{equation}
\label{eq:sub-sub}
\Diff_{2t}\supset\Diff_{[t]}\supset \Diff_{2t+\epsilon}.
\end{equation}

\sm 

{\bf \punct Final remarks on groups related to graded nilpotent Lie
algebras and
 completions of universal enveloping algebras.%
\label{ss:last}}

 1) For any finite-dimensional Lie algebra $\frh$ there is 
 a natural completion of $U(\frh)$ proposed by  Rashevsky \cite{Rash}, 1966. This completion does not contain exponentials of elements of
 $\frh$.   But it enjoys the following strong property: 
 for any unitary representation $\rho$ of the corresponding
 Lie group $H$, the Rashevsky completion 
   acts in spaces of analytic vectors of $\rho$.
This was proved by  Goodman \cite{Goo1}--\cite{Goo2} who also constructed some other completions of $U(\frh)$.

\sm 

2) Let $\frg$ be as in Subsect.~\ref{ss:1.1}. Let $\omega_1$,
\dots, $\omega_\alpha$ be a  {\it minimal} collection of homogeneous generators%
\footnote{Notice that for examples in Subsect. \ref{ss:1.1} we usually
have a natural choice of generators.} of $\frg$.
Fix $t_1$, \dots, $t_\alpha>0$. Consider the free associative algebra
$\cA_\alpha$
generated by $\omega_j$. For a monomial $\omega^I$
having degrees $d_j$ with respect to $\omega_j$
we assign the norm $\prod t_j^{d_j}$ and set that a norm
of a polynomial is the sum of norms of monomials.

Repeating the construction of Subsect. \ref{ss:sub}--\ref{ss:completions}, we get a completion $\bfU_{[t_1,\dots,t_\alpha]}(\frg)$ of $U(\frg)$
and the corresponding Banach--Lie group $\bfG_{[t_1,\dots,t_\alpha]}$.
We also can take the topological  algebra 
$$\bfU^\circ(\frg)=\bigcap_{t_1,\dots,t_\alpha}\bfU_{[t_1,\dots,t_\alpha]}(\frg)$$
and the group
 $\bfG^\circ=\cap \bfG_{[t_1,\dots,t_\alpha]}$,
 these objects do not depend on the choice of generators.

 Our considerations of Subsect. \ref{ss:sub}--\ref{ss:representations}
 literally survive in this generality (and $\bfU_{[t]}(\vvect)$
 in the notation of the previous paragraph is $\bfU_{[t,t^2]}(\vvect)$).
 
 \sm

3) The completion $\bfU_{[t_1,\dots,t_\alpha]}(\frg)$ can be characterized
in abstract terms as follows.
 For any representation $\rho$ of $\frg$ by bounded operators in a Banach space such that 
\begin{equation} 
 \tri\rho(\omega_1)\tri= t_1,\dots,  \tri\rho(\omega_\alpha)\tri= t_\alpha,
 \label{eq:ineq}
 \end{equation}
  we consider the norm
 $\tri\rho(u)\tri$ on $U(\frg)$. Next, we define
 $$
 Q_{[t_1,\dots,t_\alpha]}(u)=\sup_\rho \tri\rho(u)\tri,
 $$
 where the supremum is taken over all representation of 
 $\frg$ satisfying \eqref{eq:ineq}. We define 
  the completion of $U(\frg)$ by this norm. It coincides with the completion
 $\bfU_{[t_1,\dots,t_\alpha]}(\frg)$%
 \footnote{Indeed, $U(\frg)$ acts on the Banach space
  $\bfU_{[t_1,\dots,t_\alpha]}(\frg)$
 by right multiplications, the norms of generators are precisely
 $t_1$, \dots, $t_\alpha$, and our norm is the maximal submultiplicative 
 norm on $U(\frg)$ satisfying this property.}.
 
 We also get the topological algebra $\bfU^\circ(\frg)$
 as an intersection of all $\bfU_{[t_1,\dots,t_\alpha]}(\frg)$.
 
 \sm
 
 4)  The abstract definition of the topological algebra $\bfU^\circ(\cdot)$ 
 makes sense for a finite-dimensional Lie algebra $\frh$:
 we complete $U(\frh)$ by the family of seminorms
 $p_\rho(u)=\tri\rho(u)\tri$, where $\rho$ ranges in 
 representations of $\frh$ in Banach spaces by bounded operators.
 It was proposed by Taylor \cite{Tay}, 1972, see also
 \cite{Ari} about nilpotent $\frh$, and further references in the latter paper.
 Emphasize that our considerations above  do not survive for general finite-dimensional Lie algebras. Notice also, that the construction of $\bfU^\circ(\cdot)$
 is about representations of Lie algebras by bounded operators,
and this a priori restricts possibilities of usage of such objects (see our remarks
 in Subsect. \ref{ss:1.2}).

\tt 

University of Graz,
\\
\phantom{.}
\hfill Department of Mathematics and Scientific computing;

High School of Modern Mathematics MIPT,
\\
\phantom{.}
\hfill 1 Klimentovskiy per., Moscow; 

Moscow State University, MechMath. Dept;

 University of Vienna, Faculty of Mathematics.
 
 \sm

e-mail:yurii.neretin(dog)univie.ac.at

URL: https://www.mat.univie.ac.at/$\sim$neretin/ 

\phantom{URL:} https://imsc.uni-graz.at/neretin/index.html

\end{document}